\documentclass[12pt,a4paper]{article}

\usepackage{graphicx}
\usepackage{amsmath,amsfonts}
\usepackage{mathrsfs}
\usepackage[T1]{fontenc}
\usepackage[utf8]{inputenc}
\usepackage{authblk}
\usepackage{enumerate}
\usepackage{capt-of}
\usepackage{calc}
\usepackage{float}
\allowdisplaybreaks[1]

\newtheorem{theorem}{Theorem}
\newtheorem{lemma}{Lemma}

\let\scr\mathscr

\def\Ex{\mathbf{E}}

\def\AA{\mathbb{A}}
\def\BB{\mathbb{B}}
\def\CC{\mathbb{C}}

\def\KK{\mathbb{K}}

\def\UU{\mathbb{U}}

\def\1{\mbox{1\hspace{-.25em}I}}
\newcommand{\Liminf}{\mathop{\underline{\lim}}\limits}

\begin{document}
\title{Poisson Source Localization on the Plane. Cusp Case.}
\author[1]{O.V. Chernoyarov}
\author[2]{S. Dachian}
\author[3]{Yu.A. Kutoyants}

\affil[1,2,3]{\small  National Research University ``MPEI'', Moscow, Russia}
\affil[2]{\small University of Lille, Lille,  France}
\affil[3]{\small  Le Mans University,  Le Mans,  France}
\affil[1,2,3]{\small  Tomsk State University, Tomsk, Russia}

\date{}
\maketitle

\begin{abstract}
 This work is devoted to the problem of estimation of the localization of
Poisson source. The observations are inhomogeneous Poisson processes
 registered by the $k\geq 3$ detectors on the plane. We study the asymptotic
 properties of the  Bayes estimators in the asymptotic of large
intensities. It is supposed that the intensity functions of the  signals
arriving in the detectors  have cusp-type singularity. We show the
consistency, limit distributions and the convergence of moments of these
estimators. 
\end{abstract}

\bigskip{} \textbf{Key words}: Inhomogeneous Poisson process, Poisson
source, sensors, Bayes estimators,
cusp-type singularity.
 \bigskip{}
\date{}

\section{Introduction}

Suppose that we have $k\geq 3$ detectors at the points $D_j,j=1,\ldots,k$ with
the coordinates $\vartheta _j=\left(x_j,y_j\right),j=1,\ldots,k$ on the plane
and a source of  emission of Poisson signals at the point $D_0$ with coordinates
$\vartheta _0=\left(x_0,y_0\right)$. We consider the problem of estimation of
the position $\vartheta_0 =\left(x_0,y_0\right)$ by the observations $X
=\left(X_1,\ldots,X_k\right)$ of Poisson signals
$X_j=\left(X_j\left(t\right), 0\leq t\leq T\right)$ received by detectors
\cite{Kn10}.

An example of such model is given on the Fig. 1.


\bigskip

\begin{figure}[ht]
\hspace{3cm}\includegraphics[width=6cm]   {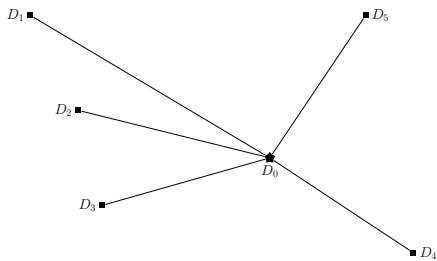}
\caption{Model of observations}
\end{figure}

 This is our third work devoted to this problem of identification
of localization of the source (see the Introduction in the work \cite{FKT18}
where we give the review of the engineering literature on this subject).

The intensity function $\lambda _{j,n}\left(\vartheta_0
,t\right)$  of the Poisson process received by the $j$-th detector taken in
this work and  in \cite{CK18}, \cite{FKT18}   is of the form
\begin{align}
\label{1}
\lambda _{j,n}\left(\vartheta_0,t\right)=n\lambda _{j}\left(t-\tau
_j\right)\1_{\left\{t\geq \tau _j\right\}}+n\lambda _0,\quad  0\leq t\leq T .
\end{align}
Here $\tau _j=\tau _j\left(\vartheta _0\right)$ is the instant of arriving of
the Poisson signal at the $j$-th detector, which is calculated by the
formula $\tau _j\left(\vartheta _0\right)=\nu ^{-1}\left\|\vartheta
_j-\vartheta _0\right\|$, where $\nu >0$ is the known rate of propagation of
the signal and $\left\|\cdot \right\|$ is Euclidean norm on the plane. The
Poisson signals are received in the presence of Poisson noise of the
known intensity $n\lambda _0>0$. The exact calculation of the error of
estimation $\Ex_{\vartheta _0}\| \bar\vartheta -\vartheta _0\|^2$ ($\bar
\vartheta $ is some estimator) in this essentially non linear statistical
problem is very difficult problem. Moreover the most interesting are the
situations where the errors of estimation are small. To obtain small errors
and have possibility to calculate it we have to consider one or another type
of asymptotics. That is why we introduce the {\it large parameter} $n$ in the
intensity function \eqref{1} and study the errors of estimation in the
asymptotics $n\rightarrow \infty $. This means that the signal and noise are
sufficiently large and the estimators $\bar\vartheta =\bar\vartheta _n$ take
values not too far from the true value: $\Ex_{\vartheta _0}\| \bar\vartheta
-\vartheta _0\|^2=o\left(1\right)$.
Recall that the similar mathematical model can be used in the problem of
GPS-localization on the plane. In this case we have $k$ emitters of the
Poisson signals and an object which receives these signals. The positions of
the emitters are known and the problem is in the estimation of the position of
the object by the observations of the signals. The intensity functions of the
received Poisson signals depend on the distance between the emitters and the
object and the receiver has to defined its position by these observations
(see, e.g. \cite{XL13}). 

 The goal of the works \cite{CK18},
\cite{FKT18} and of this one is to evaluate the errors $\Ex_{\vartheta _0}\|
\hat\vartheta -\vartheta _0\|^2 $ and $\Ex_{\vartheta _0}\| \tilde\vartheta
-\vartheta _0\|^2 $, where $\hat\vartheta _n$ is the maximum likelihood
estimator (MLE) and $\tilde\vartheta _n$ is the Bayes estimator (BE) with the
quadratic loss function. The difference between these three works is in the
conditions of regularity of the functions $\lambda _j\left(\cdot \right)$ and
as a consequence of it the rates of convergence of the errors are
different. Let us remind this class of models and errors of estimation with the
help of the Poisson process with intensity function
\begin{align*}
\lambda _{n}\left(\vartheta ,t\right)=2n\left|\frac{t-\vartheta}{\delta }
\right|^\kappa  \1_{\left\{0\leq t-\vartheta \leq \delta 
 \right\}}+2n\1_{\left\{ t \geq  \vartheta+ \delta 
 \right\}}+n,\quad 0\leq t\leq T .
\end{align*}
Here the unknown parameter $\vartheta $ is one-dimensional, $\vartheta \in \left(\alpha ,\beta \right)\subset \left[0,T\right]$.
Choosing the different values of $\kappa $ we obtain statistical problems of
different regularity. The examples of such intensities are given on the
Fig. 2, where we put $n=1$.


\bigskip

\begin{figure}[ht]
\hspace{1cm}\includegraphics[width=12cm]   {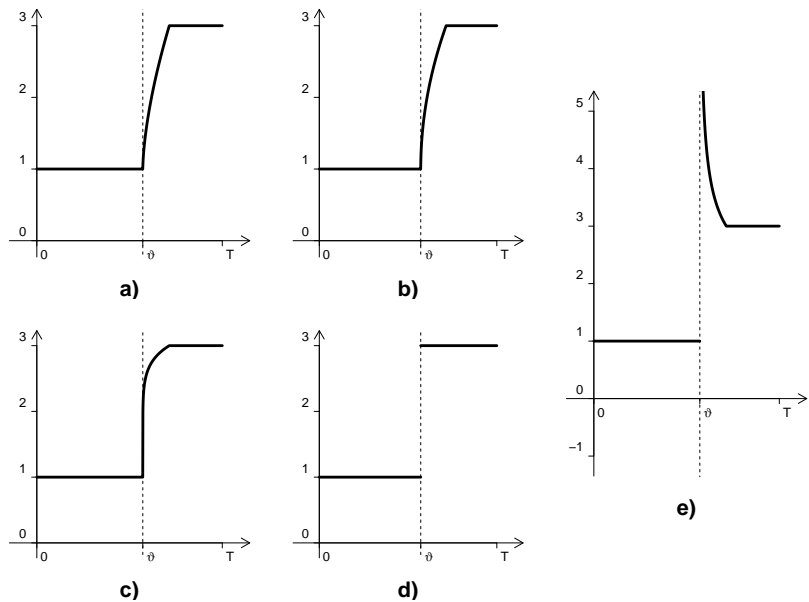}
\caption{Intensity functions of different regularity: {\bf a}) $\kappa =\frac{5}{8}$, {\bf b}) $\kappa
  =\frac{1}{2}$, {\bf c}) $\kappa =\frac{1}{8}$, {\bf d}) $\kappa =0$, {\bf e})
  $ \kappa =-\frac{3}{8}$. }
\end{figure}

The cases {\bf a}) and {\bf b}) correspond to the {\it regular} (smooth, LAN)
case. In the case {\bf c}) we have {\it cusp}-type singularity. The case
{\bf d}) corresponds to {\it change-point} model of observations and the case
{\bf e}) is {\it explosion}-type singularity. The rates of convergence of
errors in these cases are
\begin{align*}
&{\bf a})\; \Ex_{\vartheta _0}\| \tilde\vartheta _n-\vartheta _0
\|^2\approx \frac{C}{n},  &\qquad  {\bf b})\;
\Ex_{\vartheta _0}\| \tilde\vartheta _n-\vartheta _0 
\|^2\approx \frac{C}{n\ln n},\\
&{\bf c}) \; \Ex_{\vartheta _0}\| \tilde\vartheta _n-\vartheta _0
\|^2\approx \frac{C}{n^{\frac{2}{2\kappa +1}}}, &{\bf d}) \; \Ex_{\vartheta
  _0}\| \tilde\vartheta _n-\vartheta _0 \|^2\approx \frac{C}{n^{2}},\\
& {\bf e})\;\Ex_{\vartheta _0}\| \tilde\vartheta _n-\vartheta _0
\|^2\approx \frac{C}{n^{\frac{2}{2\kappa +1}}}.
\end{align*}
For the case  {\bf a}) see \cite{Kut79}, the case  {\bf b}) was considered in
\cite{CK18}, for  
the case {\bf c}) see \cite{D03}, for the case {\bf d}) see \cite{Kut84} and
the case {\bf e})  was studied in \cite{D11}.

We have to note that the study of MLE and BE in all these cases was done with
the help of some general results concerning the behavior of estimators
developed by Ibragimov and Khasminskii \cite{IH81}. Their method is based on
the study of the normalized likelihood ratio random fields, which we remind
below in this section.

   We have $k$ independent observations of 
inhomogeneous Poisson processes $X^n=\left(X_1,\ldots,X_k\right)$ with
intensities \eqref{1} depending on $\tau _j\left(\vartheta _0\right)$. We suppose that
the position of the source $\vartheta _0\in\Theta $ is unknown and we have to
estimate $\vartheta _0$ by the observations
$X^n$. Here $\Theta \subset {\cal R}^2$ is a convex
bounded and open set.

It seems that  the mathematical study of this class of models was
not yet sufficiently developed.  The statistical models of inhomogeneous
Poisson process with intensity functions having discontinuities along some
curves depending on unknown parameters were considered in \cite{Kut98},
Sections 5.2 and 5.3. Statistical inference for point processes can be found
in the works \cite{Ka91}, \cite{SM91} and \cite{St10}.

Let us recall the definitions of the MLE and BE. The functions $\lambda
_j\left(\cdot \right)$ are bounded  and the constant $\lambda >0$ therefore
the measures induced by the the processes $X_j$ in the space of their
realizations are equivalent \cite{LS01}.
The likelihood ratio function $L\left(\vartheta ,X^n\right) $  is
\begin{align*}
\ln L\left(\vartheta ,X^n\right)&=\sum_{j=1}^{k}\int_{\tau _j}^{T}\ln
\left(1+\frac{\lambda_j 
  \left(t-\tau _j\right)}{\lambda _0}\right) {\rm
  d}X_j\left(t\right) -n\sum_{j=1}^{k}\int_{\tau _j}^{T}\lambda_j
  \left(t-\tau _j\right){\rm d}t. 
\end{align*}
Of course, $\tau _j=\tau _j\left(\vartheta \right)$ and the observations
$X^n=\left(X_1,\ldots,X_n\right)$, where $X_j^n=\left(X_j\left(t\right),0\leq
t\leq T\right), j=1, \ldots,k$ are counting processes from $k$ detectors.  The
maximum likelihood estimator (MLE) $\hat\vartheta _n$ and Bayesian estimator
(BE) $\tilde\vartheta _n$ for the quadratic loss function are defined by the
``usual'' relations
\begin{align}
\label{mle}
L\left(\hat\vartheta _n ,X^n\right)=\sup_{\vartheta \in \Theta
}L\left(\vartheta ,X^n\right)
\end{align}
and
\begin{align}
\label{be}
\tilde\vartheta _n=\frac{\int_{\Theta }^{}\vartheta  p\left(\vartheta
  \right)L\left(\vartheta ,X^n\right){\rm d}\vartheta }{\int_{\Theta }^{}
  p\left(\vartheta 
  \right)L\left(\vartheta ,X^n\right){\rm d}\vartheta }.
\end{align}
Here $p\left(\vartheta \right),\vartheta \in \Theta $ is the prior
density. We suppose that it is  positive, continuous function on $\Theta
$. In this work we study the BE only. The case of MLE for this model of
observations (two-dimensional cusp) will be considered later.

\section{Main result}

Suppose that there exists a  source of Poisson signals at some point $\vartheta
_0=\left(x_0,y_0\right)\in \Theta \subset {\cal R}^2$ and $k\geq 3$ sensors
(detectors) on the same plane located at the points $\vartheta
_j=\left(x_j,y_j\right), j=1,\ldots,k$. The source was activated at the
(known) instant $t=0$ and the  signals from the source (inhomogeneous Poisson
processes) are registered by all  $k$  detectors.  The signal arrives at the
$j$-th detector at the instant $ \tau 
_j$. Of course, $\tau _j=\tau _j\left(\vartheta _0\right)$ is the  time
necessary for the signal to arrive in the $j$-th detector defined by the
relation
\begin{align*}
\tau _j\left(\vartheta _0\right)=\nu ^{-1}\left\|\vartheta _j-\vartheta
_0\right\|, 
\end{align*}
where $\nu >0$ is the known speed of propagation of the signal and
$\left\|\cdot \right\|$ is the Euclidean norm (distance) in ${\cal R}^2$.

The intensity function of the Poisson process $X_j^n=\left(X_j\left(t\right), 
0\leq t\leq T\right)$ registered by the $j$-th detector is
\begin{align}
\label{int}
\lambda _{j,n}\left(\vartheta_0 ,t\right)=nS_{j}\left(t-\tau _j\left(\vartheta _0\right)\right) +n\lambda
_0,\qquad 0\leq t\leq T, 
\end{align}
where $nS_{j}\left(t-\tau _j\left(\vartheta _0\right)\right) $ is the
intensity function of the signal and $n\lambda_0>0$ is intensity of the
noise. We suppose that the function $S_j\left(\cdot \right)$ of the signal can
be presented as follows
\begin{align}
\label{s}
S _{j}\left(t-\tau _j\right)=\lambda _j\left(t-\tau _j
\right)\left|\frac{t-\tau _j}{\delta }\right|^\kappa \1_{\left\{0\leq  t-\tau _j\leq
  \delta \right\}}+\lambda _j\left(t-\tau _j \right) \1_{\left\{t-\tau
  _j> \delta \right\}}  .
\end{align}
Here $\delta >0$ is some small parameter. This means that the signal is
strongly increasing 
function  on the interval $\left[\tau _j,\tau _j+\delta \right]$ and non
differentiable at the point $t=\tau _j$. For simplicity of the 
exposition we suppose that the noise level in all detectors is the same. 

Introduce the notations:  $\varphi _n=n^{\frac{1}{2\kappa +1}}$ and for
$j=1,\ldots,k$   
\begin{align}
\label{tau}
&\tau _j\left(\vartheta _0+\nu \varphi _n{u}\right)=\tau _j\left(\vartheta
  _0\right)- \varphi _n\langle m_j,
  u\rangle+\left\|u\right\|^2O\left(\varphi _n^2\right), \\
 & m_j=\left(\frac{x_j-x_0}{\rho _j},\frac{y_j-y_0}{\rho _j} \right)
  ,\qquad \rho _j=\left\| \vartheta _j-\vartheta _0\right\|,\quad
  \left\|m_j\right\|=1,\nonumber \\
\alpha _j&=\inf_{\vartheta \in\Theta }\tau _j\left(\vartheta \right),\quad
\beta _j=\sup_{\vartheta \in\Theta }\tau _j\left(\vartheta \right),\qquad
      {\cal T}_j=\left[\alpha _j,\beta _j\right]
.\nonumber
\end{align}

\bigskip

 Conditions   ${\scr C}$.

${\scr C}_1$. {\it The set $\Theta $ is open, convex, bounded and such that
   $0<\alpha _j<\beta _j<T$.  }

${\scr C}_2$. {\it   The source can not be in the detector, i.e.,
$\vartheta _0\not=\vartheta _j$. }

${\scr C}_3$. {\it The parameters $\kappa \in \left(0,\frac{1}{2}\right)$ and
     $\delta\in \left(0, T\right)$.}

${\scr C}_4$. {\it    The functions   
  $\lambda _j\left(t\right)>0$ have   continuous derivatives $\lambda_j
'\left(\cdot \right)$}.

${\scr C}_5$. {\it There is at least three detectors which are not on the same
  line.}

\bigskip

By the condition ${\scr C}_1$ we have $\min_j\rho _j>0$. This condition is
quite restrictive because if we take as $\Theta $ the region including
$\vartheta _0$ and all $\vartheta _j$ we have to suppose that there exists
$\varepsilon >0$ such that the discs $\CC_j=\left\{\vartheta_0
:\;\left\|\vartheta _j-\vartheta _0\right\|\leq \varepsilon \right\}$ are excluded
from $\Theta $, but in this case the set $\Theta $ is no more convex. Note
that it is possible to modify the proof in such a way that the consistency and
convergence to the limit distribution are uniform on compacts
$\KK\subset\Theta $ which do not include the positions of the detectors
$\vartheta _j$. Another point, when we do the re-normalization $\vartheta
=\vartheta_0+\nu \varphi _nu$ with $u\in \UU_n=\left\{u: \vartheta_0+\nu
\varphi _nu\in \Theta  \right\}$ we have to exclude the values $u$ which correspond
to  $\vartheta \in \CC_j$. To avoid such problems we extend the normalized
likelihood ratio random field to include these values $u$, but the true value
$\vartheta _0$ is always separated from $\vartheta _j$. 

Introduce the notations: $\lambda _j=\lambda _j\left(0\right)$,
\begin{align*}
\BB_j&=\left\{u:\;\langle m_j,u\rangle <0\right\},\qquad
\BB_j^c=\left\{u:\;\langle m_j,u\rangle\geq 0\right\},\quad \gamma
_j=\frac{\lambda _j}{\delta ^\kappa \sqrt{\lambda
    _0}},\\ J_j\left(u\right)&=J_{j,-}\left(u\right)\1_{\left\{u\in\BB_j\right\}}
+J_{j,+}\left(u\right)\1_{\left\{u\in\BB_j^c\right\}}
,\quad u\in {\cal R}^2,\\ J_{j,-}\left(u\right)&=\gamma _j\int_{0}^{\infty }
\left[\left|s+\langle m_j,u\rangle \right|^\kappa\1_{\left\{s>-\langle
    m_j,u\rangle \right\}}- \left|s\right|^\kappa \right] {\rm
  d}W_j\left(s\right),\\ J_{j,+}\left(u\right)&=\gamma _j\int_{-\langle
  m_j,u\rangle }^{\infty } \left[\left|s+\langle m_j,u\rangle \right|^\kappa-
  \left|s\right|^\kappa\1_{\left\{s>0 \right\}} \right] {\rm
  d}W_j\left(s\right),\\ R_j\left(u\right)&=R_{j,-}\1_{\left\{u\in\BB_j\right\}}
+R_{j,+}\1_{\left\{u\in\BB_j^c\right\}} 
,\quad u\in {\cal R}^2,\\ R_{j,-}&=\gamma _j^2\int_{0}^{\infty }
\left[\left|s-1 \right|^\kappa\1_{\left\{s>1 \right\}}- \left|s\right|^\kappa
  \right]^2 {\rm d}s,\\ R_{j,+}&=\gamma _j^2\int_{-1 }^{\infty }
\left[\left|s+1\right|^\kappa- \left|s\right|^\kappa\1_{\left\{s>0 \right\}}
  \right]^2 {\rm d}s.
\end{align*}
Here $W_j\left(\cdot \right), j=1,\ldots,k$ are independent Wiener processes. 
The limit likelihood ratio field is 
\begin{align*}
Z\left(u\right)=\exp\left\{\sum_{j=1}^{k}\left[J_j\left(u\right)-\frac{\left|\langle 
  m_j,u\rangle\right|^{2\kappa +1}}{2}
R_j\left(u\right)\right] \right\} ,\qquad u\in {\cal R}^2.
\end{align*}
Note that this is  a product of $k$ {\it independent } random fields
\begin{align*}
Z\left(u\right)=\prod_{j=1}^kZ_j\left(u\right),\qquad
Z_j\left(u\right)=\exp\left\{J_j\left(u\right)-\frac{\left|\langle
  m_j,u\rangle\right|^{2\kappa +1}}{2} R_j\left(u\right) \right\}.
\end{align*}
Introduce as well the random vector  $\tilde\zeta $, which  has the same
distribution as the limit of the normalized BE 
\begin{align*}
\tilde \zeta =\nu \frac{\int_{{\cal R}^2}^{}uZ\left(u\right){\rm d}u}{\int_{{\cal
      R}^2}^{}Z\left(u\right){\rm d}u} .
\end{align*}

Remark, that if all detectors are on the same line, then the 
consistent identification is impossible because the same signals come from the
symmetric with respect to this line  possible locations of the source.

 We
have the following minimax  lower bound on the mean square errors
of all estimators $\bar\vartheta _n$: Let the conditions ${\cal C}$ be
fulfilled then for any $\vartheta _0\in\Theta $
\begin{align*}
\lim_{\delta \rightarrow 0}\Liminf_{n\rightarrow \infty
}\sup_{\left\|\vartheta -\vartheta _0\right\|\leq \delta } n^{\frac{2}{2\kappa +1}}\Ex_\vartheta
\left\|\bar\vartheta _n-\vartheta \right\|^2 \geq \Ex_{\vartheta _0}
\left\|\zeta \right\|^2 .
\end{align*}
For the proof see, e.g., \cite{IH81}, Theorem 2.12.1.

We call the estimator $\bar\vartheta _n$ asymptotically efficient, if for all
$\vartheta _0\in\Theta $ we have the equality
\begin{align*}
\lim_{\delta \rightarrow 0}\lim_{n\rightarrow \infty
}\sup_{\left\|\vartheta -\vartheta _0\right\|\leq \delta } n^{\frac{2}{2\kappa +1}}\Ex_\vartheta
\left\|\bar\vartheta _n-\vartheta \right\|^2 = \Ex_{\vartheta _0}
\|\tilde \zeta \|^2 .
\end{align*}

\begin{theorem}
\label{T1} Let the conditions ${\cal R}$ be fulfilled then the BE
$\tilde\vartheta _n$ is uniformly consistent, converges in distribution
\begin{align*}
{n}^{\frac{1}{2\kappa +1}}\left(
\tilde \vartheta _n-\vartheta _0 \right)\Longrightarrow \tilde\zeta ,
\end{align*}
for any $p>0$
\begin{align*}
\lim_{n\rightarrow \infty }
n^{^{\frac{p}{2\kappa +1}}}\Ex_{\vartheta _0}\|\tilde\vartheta _n-\vartheta _0
\|^p =\Ex_{\vartheta _0}\left\| \zeta \right\|^p,
\end{align*}
 and  BE is asymptotically efficient.

\end{theorem}
{\bf Proof.} The properties of estimators mentioned in this theorem we verify
with the help of approach developed by Ibragimov and Khasminskii
\cite{IH81}. The similar method was used in the preceding our works
\cite{CK18} and \cite{FKT18}. For the convenience of understanding we  remind
it here once  more. Introduce the 
normalized likelihood ratio random field 
\begin{align*}
Z_n\left(u\right)=\frac{L\left(\vartheta_0+\nu \varphi _n{u},X^n\right)
}{L\left(\vartheta_0 ,X^n\right)},\qquad u\in
\UU_n=\left\{u:\;\vartheta_0+\nu\varphi _n{u}\in \Theta  \right\}
\end{align*}
where the normalizing function $\varphi _n=n^{-\frac{1}{2\kappa +1}}$.

Suppose that we already proved the   convergence 
\begin{align*}
Z_n\left(\cdot \right)\Longrightarrow Z\left(\cdot \right).
\end{align*}
Then the limit distribution of the BE can be  obtained as
follows (see \cite{IH81}). Below we change the variables $\vartheta
=\vartheta_u =\vartheta  
_0+ \nu\varphi _nu$.
\begin{align*}
\tilde\vartheta _n &=\frac{\int_{\Theta }^{}\theta p\left(\theta
  \right)L\left(\theta ,X^T\right){\rm d}\theta }{\int_{\Theta }^{}
  p\left(\theta \right)L\left(\theta ,X^T\right){\rm d}\theta}=\vartheta
_0+\nu\varphi _n
 \frac{\int_{\UU_n }^{}u p\left(\theta_u
  \right)L\left(\theta_u ,X^T\right){\rm d}u }{\int_{\UU_n }^{}
  p\left(\theta_u \right)L\left(\theta_u ,X^T\right){\rm d}u}
\\ &=\vartheta
_0+\nu\varphi _n\frac{\int_{\UU_n }^{}u p\left(\theta_u \right)Z_n
  \left(u\right){\rm d}u }{\int_{\UU_n }^{} p\left(\theta_u \right)Z_n
  \left(u\right){\rm d}u}.
\end{align*}
Hence
\begin{align*}
\varphi _n^{-1}\left(\tilde\vartheta _n-\vartheta
_0\right)=\nu\frac{\int_{\UU_n }^{}u p\left(\theta_u \right)Z_n
  \left(u\right){\rm d}u }{\int_{\UU_n }^{} p\left(\theta_u \right)Z_n
  \left(u\right){\rm d}u}\Longrightarrow \nu \frac{\int_{{\cal R}^2 }^{}u
  Z\left(u\right){\rm d}u }{\int_{{\cal R}^2 }^{} Z\left(u\right){\rm
    d}u}=\tilde \zeta .
\end{align*}
Recall that $p\left(\theta _u\right)\rightarrow p\left(\vartheta
_0\right)>0$.

The properties of the $Z_n\left(u\right)$ required in the Theorem
1.10.2 \cite{IH81} are checked in the three lemmas below.  Remind that this
approach to the study of the properties of these estimators was applied  in
\cite{Kut79},  \cite{Kut98}. Here we use some obtained there inequalities.

\begin{lemma}
\label{L1} Let the conditions ${\scr C}$ be fulfilled,  then the finite
dimensional distributions of the random field $Z_n\left(u\right), u\in\UU_n$
converge to the finite dimensional distributions of the limit random field
$Z\left(u\right),u\in{\cal R}^2$ and this convergence is uniform on compacts
$\KK\in\Theta $.
\end{lemma}
{\bf Proof.}  Let us denote ${\rm d}\pi _{j,n}\left(t\right)= {\rm
  d}X_j\left(t\right) -n\left[S_j\left(t-\tau _j\left(\vartheta
  _0\right)\right) +\lambda _0\right]{\rm d}t$ and put $\vartheta _u=\vartheta
_0+\nu \varphi _nu$, $\tau _j=\tau _j\left(\vartheta _0\right)$. Then we can
write
\begin{align*}
\ln Z_n\left(u\right)&=\sum_{j=1}^{k}\int_{0}^{T}\ln\left(\frac{
  S_j\left(t-\tau _j\left(\vartheta _u\right)\right)+\lambda _0
}{S_j\left(t-\tau _j\right)+\lambda _0
}\right){\rm d}\pi _{j,n}\left(t\right)\\
&\qquad -n\sum_{j=1}^{k}\int_{0}^{T}\left[\frac{S
  _j\left(t-\tau _j\left(\vartheta _u\right)\right)+\lambda _0}{
      S_j\left(t-\tau _j\right)+\lambda _0 }  
  -1\right.\\
&\qquad \left.- \ln\left(\frac{
  S_j\left(t-\tau _j\left(\vartheta _u\right)\right)+\lambda _0
}{S_j\left(t-\tau _j\right)+\lambda _0
}\right)\right] \left[S_j\left(t-\tau _j\right)+\lambda _0 \right]{\rm d}t\\
&=\sum_{j=1}^{k}\int_{0}^{T}F_j\left(t,\vartheta _u\right){\rm d}\pi
_{j,n}\left(t\right)-n\sum_{j=1}^{k}\int_{0}^{T}G_j\left(t,\vartheta _u\right){\rm d}t 
\end{align*}
with obvious notation.

Let $u\in {\BB}_j$. Then $\tau _j\left(\vartheta _u\right)>\tau _j $. Folowing
the same arguments as that given in \cite{D03}, we obtain the 
asymptotic ($n\rightarrow \infty $) relations: 
\begin{align*}
J_{j,n}\left(u\right)&=\int_{0}^{T}F_j\left(t,\vartheta _u\right){\rm d}\pi
_{j,n}\left(t\right)=\int_{\tau _j}^{\tau
    _j+\delta }F_j\left(t,\vartheta _u\right){\rm
    d}\pi_{j,n}\left(t\right) \left(1+o\left(1\right)\right)
 \\
I_{j,n}\left(u\right)&=n\int_{0}^{T}G_j\left(t,\vartheta _u\right){\rm d}t=n\int_{\tau _j}^{\tau
    _j+\delta }G_j\left(t,\vartheta _u\right){\rm
    d}t \left(1+o\left(1\right)\right)
\end{align*}

For  $t\in \left[\tau _j,\tau _j-\varphi _n\langle m_j,u\rangle \right]$ as $\varphi _n\rightarrow 0$ we obtain the expansions
\begin{align*}
\lambda _j(t-\tau
  _j(\vartheta _u))&=\lambda _j(0)+\left(t-\tau
  _j(\vartheta _u)\right) \lambda _j'(0)\left(1+o\left(1\right)\right)=\lambda _j+o\left(1\right) \\
\lambda _j(t-\tau
  _j(\vartheta _u))&=\lambda _j(t-\tau
  _j) +\varphi _n \langle m_j,u\rangle \lambda _j'(t-\tau
  _j)+O\left(\varphi _n^2\right)\left\|u\right\|^2, \\
\left| \frac{t-\tau _j\left(\vartheta _u\right)}{\delta }\right|^\kappa
&=\delta ^{-\kappa } \left|t-\tau _j+\varphi
_n\langle m_j,u\rangle +O\left(\varphi _n^2\right)\right|^\kappa \\
&=\delta ^{-\kappa } \left|t-\tau _j+\varphi
_n\langle m_j,u\rangle \right|^\kappa +O\left(\varphi _n^{2\kappa }\right)
\end{align*}
Here we used the inequality $\left|a+b\right|^\kappa \leq
\left|a\right|^\kappa+\left|b\right|^\kappa$.

 Further, for $\tau _j\leq t\leq \tau _j-\varphi _n\langle m_j,u\rangle $ and 
$\left\|u\right\|<L$ we can write
\begin{align*}
&\ln\left(\frac{S_j\left(t-\tau _j\left(\vartheta _u\right))+\lambda _0\right)
  }{S_j\left(t-\tau _j\right)+\lambda
    _0}\right)=\ln\left(\frac{\lambda _0 }{\lambda _j\left(t-\tau
    _j\right)\left|\frac{t-\tau _j}{\delta }\right|^\kappa +\lambda _0}\right)\\
 &\qquad \qquad =-\ln\left(1+\frac{\lambda _j}{\lambda _0}\left|\frac{t-\tau
    _j}{\delta }\right|^\kappa \right)\left(1+O\left(\varphi _n\right)\right)\\ 
 &\qquad \qquad =-\frac{\lambda _j}{\lambda _0}\left|\frac{t-\tau _j}{\delta
  }\right|^\kappa \left(1+O\left(\varphi _n^{2\kappa }\right)\right).
\end{align*}
For $t\in \left[\tau _j-\varphi _n\langle m_j,u\rangle ,\delta  \right]$ the
similar relations are 
\begin{align*}
&\ln\left(\frac{S_j\left(t-\tau _j\left(\vartheta _u\right))+\lambda _0\right) 
  }{S_j\left(t-\tau _j\right)+\lambda
    _0}\right)=\ln\left(\frac{  \lambda _j\left(t-\tau
    _j\left(\vartheta _u\right)\right)\left|\frac{t-\tau _j\left(\vartheta
      _u\right)}{\delta }\right|^\kappa +     \lambda _0 }{\lambda
    _j\left(t-\tau 
    _j\right)\left|\frac{t-\tau _j}{\delta }\right|^\kappa +\lambda
    _0}\right)\\ 
 &\quad =\ln\left(1+\frac{\lambda _j\left(t-\tau _j\left(\vartheta _u\right)
    \right)\left|\frac{t-\tau 
    _j\left(\vartheta _u\right)}{\delta }\right|^\kappa-\lambda
    _j\left(t-\tau _j \right)\left|\frac{t-\tau 
    _j}{\delta }\right|^\kappa}{\lambda
    _j\left(t-\tau _j \right)\left|\frac{t-\tau 
    _j}{\delta }\right|^\kappa+\lambda _0}\right)\\ 
 &\qquad \qquad =\frac{\lambda _j}{\lambda _0}\left[\left|\frac{t-\tau
      _j+\varphi _n\langle m_j,u\rangle }{\delta
  }\right|^\kappa- \left|\frac{t-\tau _j}{\delta
  }\right|^\kappa    \right]\left(1+O\left(\varphi _n^{2\kappa
  }\right)\right). 
\end{align*}
Therefore
\begin{align*}
&\Ex_{\vartheta _0}
  \left(J_{j,n}\left(u\right)\right)^2=\int_{0}^{T}F_j\left(t,\vartheta
  _u\right)^2\lambda 
  _{j,n}\left(\vartheta _0,t\right){\rm d}t\\
&\qquad =\frac{\lambda _j^2n}{\lambda _0}   \int_{\tau _j}^{\tau _j-\varphi
    _n\langle m_j,u\rangle } \left|\frac{t-\tau _j}{\delta
  }\right|^{2\kappa}{\rm d}t+o\left(1\right)\\
&\qquad \quad +\frac{\lambda _j^2n}{\lambda _0}   \int_{\tau _j-\varphi
    _n\langle m_j,u\rangle}^{\delta }\left[\left|\frac{t-\tau
      _j+\varphi _n\langle m_j,u\rangle }{\delta
  }\right|^\kappa- \left|\frac{t-\tau _j}{\delta
  }\right|^\kappa    \right]^2{\rm d}t\\
&\qquad =\frac{\lambda _j^2n}{\lambda _0\delta ^{2\kappa }} \varphi
  _n^{2\kappa +1}  \int_{0}^{-\langle m_j,u\rangle }
  \left|s\right|^{2\kappa}{\rm d}s\\ 
&\qquad \quad +\frac{\lambda _j^2n}{\lambda _0\delta ^{2\kappa }}\varphi
  _n^{2\kappa +1}  
  \int_{-\langle m_j,u\rangle}^{\frac{\delta-\tau _j}{\varphi _n}
  }\left[\left|s+\langle m_j,u\rangle \right|^\kappa- \left|s\right|^\kappa 
    \right]^2{\rm d}s+o\left(1\right)\\
&\qquad =\gamma _j^2  \left|\langle m_j,u\rangle \right|
  ^{2\kappa +1}  \int_{0}^{1}
  \left|v\right|^{2\kappa}{\rm d}v\\ 
&\qquad \quad +\gamma _j^2 \left|\langle m_j,u\rangle \right|
  ^{2\kappa +1}
  \int_{1}^{-\frac{\delta-\tau _j}{ \langle m_j,u\rangle  \varphi _n}
  }\left[\left|v-1\right|^\kappa- \left|v\right|^\kappa 
    \right]^2{\rm d}v+o\left(1\right)\\ 
&\qquad \quad =\gamma _j^2 \left|\langle m_j,u\rangle \right|
  ^{2\kappa +1}
  \int_{0}^{-\frac{\delta-\tau _j}{ \langle m_j,u\rangle  \varphi _n}
  }\left[\left|v-1\right|^\kappa\1_{\left\{v\geq 1\right\}}- \left|v\right|^\kappa 
    \right]^2{\rm d}v+o\left(1\right)\\ 
&\qquad \quad =\gamma _j^2 \left|\langle m_j,u\rangle \right|
  ^{2\kappa +1}R_n+o\left(1\right),
\end{align*}
where we changed the variables $t=\tau _j+s\varphi _n$ and  $s=-v\langle m_j,u\rangle $.
Recall that $n\varphi _n^{2\kappa +1}=\nu^{2\kappa +1} $ and 
$\gamma _j^2={\lambda _j^2\nu ^{2\kappa +1}}{\lambda _0^{-1}\delta ^{-2\kappa }}
$. Hence  for $u\in\BB_-$  we obtain the following limit
\begin{align*}
R_n&=\int_{0}^{-\frac{\delta-\tau _j}{ \langle m_j,u\rangle  \varphi _n}
  }\left[\left|v-1\right|^\kappa\1_{\left\{v\geq 1\right\}}- \left|v\right|^\kappa 
    \right]^2{\rm d}v \\
&\qquad \longrightarrow \int_{0}^{\infty 
  }\left[\left|v-1\right|^\kappa\1_{\left\{v\geq 1\right\}}- \left|v\right|^\kappa 
    \right]^2{\rm d}v =R_{j,-}.
\end{align*}
These arguments allow us to write the representation
\begin{align*}
J_{j,n}\left(u\right)&=\gamma _j\int_{0}^{\frac{\delta-\tau _j}{\varphi _n} }
\left[\left|s+\langle
  m_j,u\rangle \right|^\kappa\1_{\left\{s\geq
    -\langle m_j,u\rangle \right\}}- \left|s\right|^\kappa \right]{\rm 
  d}W_{j,n}\left(s\right)+o\left(1\right).
\end{align*}
Here
\begin{align*}
&W_{j,n}\left(s\right)=\frac{1}{\sqrt{\lambda _0n\varphi _n}} \left[X_j\left(\tau
    _j+s\varphi 
  _n\right)-X_j\left(\tau _j\right)-\int_{\tau _j}^{\tau _j+s\varphi
    _n}\lambda _{j,n}\left( \vartheta _0,v\right){\rm d}v  \right],\\
&\Ex_{\vartheta
  _0}W_{j,n}\left(s\right)^2= \frac{n}{\lambda _0n\varphi _n}\int_{\tau _j}^{\tau
    _j+s\varphi     _n}\lambda _{j,n}\left( \vartheta _0,v\right){\rm
    d}v=s+o\left(1\right),\\ 
&\Ex_{\vartheta _0}W_{j,n}\left(s\right)=0,\quad \quad \Ex_{\vartheta
    _0}W_{j,n}\left(s_1\right)W_{j,n}\left(s_2\right)=s_1\wedge s_2+o\left(1\right) 
\end{align*}

The standard central limit theorem  provides us the corresponding convergence
of stochastic integrals. For any $u_1,\ldots, u_M\in \BB_j$ we have the joint
asymptotic normality of the stochastic integrals
\begin{align*}
Y_{j,n}\equiv \Bigl(J_{j,n}\left(u_1\right),\ldots,
J_{j,n}\left(u_M\right)\Bigr)\Longrightarrow Y_{j}\equiv
\Bigl(J_{j}\left(u_1\right),\ldots, J_{j}\left(u_M\right)\Bigr),
\end{align*}
where
\begin{align*}
J_j\left(u\right)&=\gamma _j\int_{0}^{\infty  }
\left[\left|s+\langle
  m_j,u\rangle \right|^\kappa\1_{\left\{s\geq
    -\langle m_j,u\rangle \right\}}- \left|s\right|^\kappa \right] {\rm 
  d}W_j\left(s\right).
\end{align*}
Moreover, the similar arguments gives us the convergence
\begin{align}
\label{k-m}
{\bf Y}_n\equiv \Bigl(Y_{1,n},\ldots,Y_{k,n}\Bigr)\Longrightarrow {\bf Y}\equiv \Bigl(Y_{1},\ldots,Y_{k}\Bigr)
\end{align}

Consider now the values $u\in \BB_j^c$. Then $\tau _j\left(\vartheta
_u\right)\leq \tau _j\left(\vartheta _0\right)$ or asymptotically $\tau
_j\left(\vartheta _0\right) -\varphi _n\langle m_j,u\rangle +O\left(\varphi
_n^2\right)\leq \tau _j\left(\vartheta _0\right)$. The similar arguments allow
us to verify the convergence \eqref{k-m} with the limit process
\begin{align*}
J_j\left(u\right)&=\gamma _j\int_{-\langle
  m_j,u\rangle }^{\infty  }\left[\left|s+\langle
  m_j,u\rangle \right|^{\kappa}\1_{\left\{s\leq 0
    \right\}}  +\left[\left|s+\langle
  m_j,u\rangle \right|^\kappa- \left|s\right|^\kappa \right]\right] {\rm 
  d}W\left(s\right).
\end{align*}

Therefore we have the convergence of finite-dimensional distributions of the
stochastic integrals.

For the ordinary integral $ I_{j,n}\left(u\right)$ we have the similar
representation ($u\in \BB, G_{j,t}=G_j\left(t,u\right)$)
\begin{align*}
I_{j,n}\left(u\right)&=n\int_{0}^{\tau _j}G_{j,t}\;{\rm d}t+n\int_{\tau
  _j}^{\tau _j\left(\vartheta _u\right)}G_{j,t}\;{\rm d}t+n\int_{\tau
  _j\left(\vartheta _u\right)}^{\tau _j+\delta
}G_{j,t}\;{\rm d}t\\
&\qquad \qquad +n \int_{\tau _j+\delta}^{\tau _j\left(\vartheta
  _u\right)+\delta}G_{j,t}\;{\rm d}t+n \int_{\tau _j\left(\vartheta
  _u\right)+\delta}^{T }G_{j,t}\;{\rm d}t\\
&\quad =n\int_{\tau
  _j}^{\tau _j\left(\vartheta _u\right)}G_{j,t}\;{\rm d}t +n\int_{\tau
  _j\left(\vartheta _u\right)}^{\tau _j+\delta
}G_{j,t}\;{\rm d}t+n\int_{\tau
  _j+\delta }^{\tau _j\left(\vartheta _u\right)+\delta
}G_{j,t}\;{\rm d}t+o\left(1\right)
\end{align*}
For $t\in \left[0,\tau _j \right]$ we have $ G_j\left(t,u\right)=0$ and for
$t\in \left[ \tau _j+\delta ,T\right]$ the function
$G_j\left(t,u\right) $ has continuous bounded derivative and we can write
\begin{align*}
n \int_{\tau _j\left(\vartheta
  _u\right)+\delta}^{T }G_{j,t}\;{\rm d}t\leq Cn\varphi _n^2\left\|u\right\|^2=o\left(1\right).
\end{align*}
 Consider the case $t\in \left[\tau _j,\tau _j\left(\vartheta _u\right)
  \right]$. Using expansion
$\ln\left(1+x\right)=x-\frac{x^2}{2}+O\left(x^3\right) $ we can write
\begin{align*}
&\frac{\lambda _{j,n}\left(\vartheta _u,t\right)}{\lambda
    _{j,n}\left(\vartheta _0,t\right) } -1- \ln\left(\frac{\lambda
    _{j,n}\left(\vartheta _u,t\right)}{\lambda _{j,n}\left(\vartheta
    _0,t\right) } \right)=\frac{\lambda _0}{\lambda _j\left(t-\tau
    _j\right)\left|\frac{t-\tau _j}{\delta }\right|^\kappa +\lambda
    _0}\\ &\quad \qquad -1-\ln \left(\frac{\lambda _0}{\lambda _j\left(t-\tau
    _j\right)\left|\frac{t-\tau _j}{\delta }\right|^\kappa +\lambda _0}
  \right) = \frac{\lambda _j^2}{2\lambda _0^2\delta ^{2\kappa }}\left|{t-\tau
    _j}\right|^{2\kappa }\left(1+o\left(1\right)\right).
\end{align*}
For $t\in \left[\tau _j\left(\vartheta _u\right) ,\tau _j +\delta \right]$ we have
\begin{align*}
&\frac{\lambda _{j,n}\left(\vartheta _u,t\right)}{\lambda _{j,n}\left(\vartheta
  _0,t\right) } -1- \ln\left(\frac{\lambda _{j,n}\left(\vartheta
  _u,t\right)}{\lambda _{j,n}\left(\vartheta _0,t\right) }
\right)=\frac{\lambda _j\left(t-\tau _j\left(\vartheta
  _u\right)\right)\left|\frac{t-\tau 
    _j \left(\vartheta _u\right)}{\delta }\right|^\kappa +\lambda _0}{\lambda
  _j\left(t-\tau _j\right)\left|\frac{t-\tau 
    _j}{\delta }\right|^\kappa +\lambda _0}\\
 &\quad \qquad -1-\ln
\left(\frac{\lambda _j\left(t-\tau _j\left(\vartheta
  _u\right)\right)\left|\frac{t-\tau 
    _j \left(\vartheta _u\right)}{\delta }\right|^\kappa +\lambda _0}{\lambda
  _j\left(t-\tau _j\right)\left|\frac{t-\tau 
    _j}{\delta }\right|^\kappa +\lambda _0} \right)\\
&\qquad  = \frac{\lambda
  _j^2}{2\lambda _0^2\delta ^{2\kappa }   }\Bigr(\left|{t-\tau _j\left(\vartheta
  _u\right)}\right|^{\kappa 
}-\left|{t-\tau _j}\right|^{\kappa
}\Bigr)^2\left(1+o\left(1\right)\right).
\end{align*}
These relations allow us to write
\begin{align*}
I_{j,n}&=\frac{n\lambda
  _j^2}{2\lambda _0\delta ^{2\kappa }   }\int_{\tau
  _j}^{\tau _j\left(\vartheta _u\right)} \left|{t-\tau
    _j}\right|^{2\kappa }     {\rm d}t\\
&\qquad  +\frac{n\lambda
  _j^2}{2\lambda _0\delta ^{2\kappa }   }\int_{\tau
  _j\left(\vartheta _u\right)}^{\tau _j+\delta
}\Bigr(\left|{t-\tau _j\left(\vartheta
  _u\right)}\right|^{2\kappa 
}-\left|{t-\tau _j}\right|^{2\kappa
}\Bigr){\rm d}t+o\left(1\right)\\
&=\frac{n\varphi _n^{2\kappa +1}\lambda
  _j^2}{2\lambda _0\delta ^{2\kappa }   }\int_{0}^{-\langle m_j,u\rangle} \left|{s}\right|^{2\kappa }     {\rm d}s\\
&\qquad\qquad   +\frac{n\lambda
  _j^2\varphi _n^{2\kappa +1}}{2\lambda _0\delta ^{2\kappa }   }\int_{-\langle
  m_j,u\rangle}^{\frac{\delta }{\varphi _n} 
}\Bigr(\left|s+\langle m_j,u\rangle \right|^{\kappa 
}-\left|s\right|^{\kappa
}\Bigr)^2{\rm d}s+o\left(1\right)\\
\\
&=\frac{\lambda
  _j^2\left|\langle m_j,u\rangle\right|^{2\kappa +1}}{2\lambda _0\delta ^{2\kappa }   }\int_{0}^{1} \left|{s}\right|^{2\kappa }     {\rm d}s\\
&\qquad\qquad   +\frac{n\lambda
  _j^2\varphi _n^{2\kappa +1}}{2\lambda _0\delta ^{2\kappa }   }\int_{1}^{\tau
  _j-\tau _j\left(\vartheta _u\right)   \frac{\delta }{\varphi _n}  
}\Bigr(\left|v-1 \right|^{\kappa 
}-\left|v\right|^{\kappa
}\Bigr)^2{\rm d}s+o\left(1\right)\\
&\rightarrow \frac{\gamma _j^2}{2}\int_{0}^{\infty }\left[
  \left|s\right|^{2\kappa } \1_{\left\{s< -\langle m_j,u\rangle\right\}}+\Bigr(\left|s+\langle m_j,u\rangle \right|^{\kappa 
}-\left|s\right|^{\kappa}\Bigr)^2\1_{\left\{s\geq  -\langle m_j,u\rangle\right\}}\right]{\rm d}s
.
\end{align*}


Note that all convergences mentioned above are uniform on compacts
$\KK\subset\Theta $.

\begin{lemma} 
\label{L2}
 Let the condition ${\scr R}_2 $ be fulfilled,  then there exists a
  constant $C>0$, which does not depend on $n$ such that for any $R>0$
\begin{align*}
\sup_{\vartheta _0\in\Theta }\sup_{\left\|u_1\right\|+\left\|u_2\right\|\leq
  R}\left\|u_1-u_2\right\|^{-2\kappa-1 }\Ex_{\vartheta _0} \left|
Z_n^{\frac{1}{2}}\left(u_1\right)-Z_n^{\frac{1}{2}}\left(u_2\right)\right|^2\leq 
C \left(1+R\right) .
\end{align*}
\end{lemma}
{\bf Proof.} 
We have the estimate (see, e.g. \cite{Kut98})
\begin{align*}
&\Ex_{\vartheta _0} \left|
  Z_n^{\frac{1}{2}}\left(u_1\right)-Z_n^{\frac{1}{2}}\left(u_2\right)\right|^2\leq
\sum_{j=1}^{k}  \int_{0}^{T}\left[\sqrt{\lambda _{j,n} \left(\vartheta _{u_2}
      ,t\right)}-\sqrt{\lambda _{j,n} \left(\vartheta _{u_1},t\right)}
    \right]^2{\rm d}t\\ 
&\qquad =\sum_{j=1}^{k}\int_{0}^{T}\frac{n^2\left[S
      _{j} \left(t-\tau _j\left(\vartheta
        _{u_2}\right)\right)-{S _{j} \left(t-\tau _j\left(\vartheta
        _{u_1}\right)\right)} \right]^2}{\left[\sqrt{\lambda _{j,n} \left(\vartheta
        _{u_2} ,t\right)}+\sqrt{\lambda _{j,n} \left(\vartheta _{u_1},t\right)}
      \right]^2}{\rm d}t\\
 &\qquad \leq \frac{n}{4\lambda
    _0}\sum_{j=1}^{k}\int_{0}^{T}{\left[S
      _{j} \left(t-\tau _j\left(\vartheta
        _{u_2}\right)\right)-{S _{j} \left(t-\tau _j\left(\vartheta
        _{u_1}\right)\right)} \right]^2}{\rm d}t,
\end{align*}
where we used the estimate $ \lambda _{j,n} \left(\vartheta _{u},t\right)\geq
n\lambda _0$. Suppose that ${\tau _j( \vartheta _{u_1})}<{\tau _j( \vartheta
  _{u_2})}$ and denote $ \Delta _t =\sqrt{n}\left[S
      _{j} \left(t-\tau _j\left(\vartheta
        _{u_2}\right)\right)-{S _{j} \left(t-\tau _j\left(\vartheta
        _{u_1}\right)\right)}\right]$. Then 
\begin{align*}
\int_{0}^{T}\Delta _t^2{\rm d}t&=\int_{0}^{\tau _j( \vartheta _{u_1})}\Delta
_t^2{\rm d}t +\int_{\tau _j( \vartheta _{u_1})}^{\tau _j( \vartheta
  _{u_2})}\Delta _t^2{\rm d}t+\int_{\tau _j( \vartheta _{u_2})}^{T}\Delta _t^2
    {\rm d}t \\ & =\int_{\tau _j( \vartheta _{u_1})}^{\tau _j( \vartheta
      _{u_2})}\Delta _t^2 {\rm d}t+\int_{\tau _j( \vartheta _{u_2})}^{T}\Delta
    _t^2 {\rm d}t .
\end{align*}
Remark that  the function $\Delta _t=0$ on the interval $\left[0,{\tau _j(
    \vartheta _{u_1})}\right]$ and $\Delta _t={nS _{j} \left(\vartheta _{u_1},t\right)} $
on the interval $\left[{\tau _j(
    \vartheta _{u_1})},{\tau _j( 
    \vartheta _{u_2})}\right]$. Therefore
\begin{align*}
\int_{\tau _j( \vartheta _{u_1})}^{\tau _j(
    \vartheta _{u_2})}\Delta _t^2{\rm d}t&=n\int_{\tau _j( \vartheta _{u_1})}^{\tau _j(
    \vartheta _{u_2})}\lambda _{j}
\left(t-\tau _j(\vartheta _{u_1})\right) ^2\left|\frac{t-\tau _j(\vartheta
  _{u_1})}{\delta }\right|^{2\kappa }{\rm d}t\\ 
& \leq Cn\int_{\tau _j( \vartheta _{u_1})}^{\tau _j(
    \vartheta _{u_2})}\left|\frac{t-\tau _j(\vartheta
  _{u_1})}{\delta }\right|^{2\kappa }{\rm d}t\leq Cn\left|\frac{\tau _j(\vartheta
  _{u_2})-\tau _j(\vartheta
  _{u_1})}{\delta }\right|^{2\kappa+1 }\\
& \leq Cn\varphi _n^{2\kappa+1 }\left\|u_2-u_1\right\|^{2\kappa+1 }=
C\left\|u_2-u_1\right\|^{2\kappa+1 } .
\end{align*}
Further
\begin{align}
\label{si}
\int_{\tau _j( \vartheta _{u_2})}^{T}\Delta _t^2 {\rm d}t=\int_{\tau _j(
  \vartheta _{u_2})}^{\tau _j( \vartheta _{u_1})+\delta }\Delta _t^2 {\rm
  d}t+\int_{\tau _j( \vartheta _{u_1})+\delta }^{\tau _j( \vartheta
  _{u_2})+\delta }\Delta _t^2 {\rm d}t+\int_{\tau _j( \vartheta _{u_2})+\delta
}^{T }\Delta _t^2 {\rm d}t.
\end{align}
Using the estimate   
\begin{align*}
\left|\lambda _{j} \left(t-\tau _j(\vartheta _{u_2})\right)-\lambda _{j}
\left(t-\tau _j(\vartheta _{u_1})\right)\right|^2\leq C\varphi
_n^2\left\|u_2-u_1\right\|^2 
\end{align*}
we obtain for   the first integral
\begin{align*}
& \int_{\tau _j( \vartheta _{u_2})}^{\tau _j( \vartheta _{u_1})+\delta
}\Delta_t^2 {\rm d}t=n\int_{\tau _j( \vartheta _{u_2})}^{\tau _j( \vartheta
  _{u_1})+\delta } \left[ \lambda _{j}
\left(t-\tau _j(\vartheta _{u_2})\right) \left|\frac{t-\tau _j(\vartheta
  _{u_2})}{\delta }\right|^{\kappa }\right.\\
&\qquad \qquad \left.-\lambda _{j}
\left(t-\tau _j(\vartheta _{u_1})\right)\left|\frac{t-\tau _j(\vartheta
  _{u_1})}{\delta }\right|^{\kappa }    \right]^2{\rm d}t\\
&\qquad \leq Cn\varphi
_n^2\left\|u_2-u_1\right\|^2 +Cn\int_{\tau _j( \vartheta _{u_2})}^{\tau _j( \vartheta
  _{u_1})+\delta } \left[  \left|{t-\tau _j(\vartheta
  _{u_2})}\right|^{\kappa }-\left|{t-\tau _j(\vartheta
  _{u_1})}\right|^{\kappa }    \right]^2{\rm d}t\\
&\qquad \leq
 C\varphi_n^{1-2\kappa }\left\|u_2-u_1\right\|^2 \\
&\qquad\qquad \quad 
+Cn\varphi _n^{2\kappa +1}\int_{0}^{\frac{\tau _j( \vartheta _{u_1})-\tau _j(
    \vartheta _{u_2})+\delta}{\varphi _n} } 
 \left[  \left|s\right|^{\kappa}
      -\left|s-\frac{\tau _j(\vartheta   _{u_1})-\tau _j(\vartheta
      _{u_2})}{\varphi _n}\right|^{\kappa }    \right]^2{\rm d}s\\
&\qquad \leq
 C\varphi_n^{1-2\kappa }\left\|u_2-u_1\right\|^2+C \left\|u_2-u_1\right\|^{2\kappa +1},
\end{align*}
where we used the relations 
\begin{align*}
&\left|\frac{\tau _j(\vartheta_{u_1})-\tau _j(\vartheta_{u_2})}{\varphi _n}+\langle m_j,u_1\rangle
-\langle m_j,u_2\rangle\right|  \leq C\varphi _n \left\|u_2-u_1\right\|^2,\\
&\int_{0}^{\infty }
 \left[  \left|s\right|^{\kappa}
      -\left|s-\langle m_j,u_1-u_2\rangle
\right|^{\kappa }    \right]^2{\rm d}s\\
&\qquad \quad \leq \left\|u_2-u_1\right\|^{2\kappa +1}
 \int_{0}^{\infty } 
 \left[  \left|v\right|^{\kappa}
      -\left|v-\langle m_j,e\rangle
\right|^{\kappa }    \right]^2{\rm d}v\leq C \left\|u_2-u_1\right\|^{2\kappa +1}.
\end{align*}
Here we set $s=v\left\|u_2-u_1\right\|$ and
$e=\left\|u_2-u_1\right\|^{-1}\left( u_2-u_1\right)$. 

As the function $S\left(t\right)$
has  a bounded derivative $S'\left(t\right)$ on the interval $\left[{\tau _j(
    \vartheta _{u_2})+\delta 
},{T } \right]$ we can write

\begin{align*}
\int_{\tau _j( \vartheta _{u_2})+\delta }^{T }\Delta _t^2 {\rm d}t\leq
Cn\varphi _n^2\left\|u_2-u_1\right\|^2 \leq C\left(1+R\right)\left\|u_2-u_1\right\|^{2\kappa +1}.
\end{align*}

The other cases  can be estimated by a similar way.

\begin{lemma} 
\label{L3}
Let the conditions ${\scr C} $ be fulfilled,  then there exists a
  constant $\kappa >0$, which does not depend on $n$ such that
\begin{align}
\label{pi}
\sup_{\vartheta _0\in\Theta }\Ex_{\vartheta _0} 
Z_n^{\frac{1}{2}}\left(u\right)\leq e^{-\kappa \left\|u\right\|^{\frac{2}{2\kappa +1}}}.
\end{align}
\end{lemma}
{\bf Proof.} Let us denote $\theta _u=\vartheta _0+\nu \varphi _nu$ and put 
\begin{align*}
Z_{j,n}\left(u\right)&=\exp\left\{\int_{0}^{T}\ln\left(\frac{\lambda
  _{j,n}\left(\theta_u,t\right)}{\lambda
  _{j,n}\left(\vartheta_0,t\right)}\right){\rm
  d}X_j\left(t\right)\right.\\
&\qquad \qquad \qquad \left.-\int_{0}^{T}\left[\lambda
  _{j,n}\left(\theta_u,t\right)-\lambda
  _{j,n}\left(\vartheta_0,t\right) \right]{\rm d}t\right\} .
\end{align*} 
By  Lemma 2.2  in \cite{Kut98} we can write 
\begin{align*}
\Ex_{\vartheta _0} Z_{j,n}^{\frac{1}{2}}\left(u\right)
=\exp\left\{-\frac{1}{2}\int_{0}^{T}\left[\sqrt{\lambda
  _{j,n}\left(\theta_u,t\right)}- \sqrt{\lambda
  _{j,n}\left(\vartheta_0,t\right)} \right]^2{\rm d}t       \right\} .
\end{align*}
Hence
\begin{align}
\Ex_{\vartheta _0} Z_n^{\frac{1}{2}}\left(u\right)&=\prod_{j=1}^k
\Ex_{\vartheta _0} Z_{j,n}^{\frac{1}{2}}\left(u\right)\nonumber\\ 
&=\exp\left\{-\frac{1}{2}\sum_{j=1}^{k}\int_{0}^{T}\left[\sqrt{\lambda
  _{j,n}\left(\theta_u,t\right)}- \sqrt{\lambda
  _{j,n}\left(\vartheta_0,t\right)} \right]^2{\rm d}t       \right\} .
\label{32}
\end{align}
First for simplicity of calculation we write
\begin{align}
&\int_{0}^{T}\left[\sqrt{\lambda
  _{j,n}\left(\vartheta,t\right)}- \sqrt{\lambda
  _{j,n}\left(\vartheta_0,t\right)} \right]^2{\rm d}t \nonumber\\
&\qquad\qquad\qquad= \int_{0}^{T}\frac{\left[{\lambda
  _{j,n}\left(\vartheta,t\right)}- {\lambda
  _{j,n}\left(\vartheta_0,t\right)} \right]^2}{\left[\sqrt{\lambda
  _{j,n}\left(\vartheta,t\right)}+ \sqrt{\lambda
  _{j,n}\left(\vartheta_0,t\right)} \right]^2}{\rm d}t \nonumber\\
&\qquad\qquad\qquad\geq c_jn\int_{0}^{T}\left[{S
  _{j}\left(t-\tau _j\left(\vartheta\right)\right)}- {S
  _{j}\left(t-\tau _j\right)} \right]^2{\rm d}t ,
\label{bel}
\end{align}
where $c_j=\left(4\lambda _M\right)^{-1}>0$ and $\lambda _M=\lambda
_0+\max_{t\in {\cal T}_j}
S_j(t)$. Therefore it is sufficient to study the integral
\begin{align*}
I_j\left(\vartheta \right)&=\int_{0}^{T}\left[{S
  _{j}\left(t-\tau _j\left(\vartheta\right)\right)}- {S
  _{j}\left(t-\tau _j\right)} \right]^2{\rm
  d}t\\
&=\int_{\tau _j\left(\vartheta \right)\wedge \tau _j}^{T}\left[{S
  _{j}\left(t-\tau _j\left(\vartheta\right)\right)}- {S
  _{j}\left(t-\tau _j\right)} \right]^2{\rm d}t.
\end{align*}

We evaluate  these integrals on  two  sets $\AA=\left\{\vartheta
:\;\left\|\vartheta -\vartheta _0\right\|\leq h\right\}$ and $\AA^c$. Here
$h>0$ is some small number. Recall that we denoted $\tau _j=\tau
_j\left(\vartheta _0\right)$. 

Let $\vartheta \in \AA\cap \BB$, where $ \BB=\left\{\vartheta \in \AA:\;
\tau_j \left(\vartheta \right)>\tau _j\left(\vartheta _0\right) \right\}
$. Moreover $\tau_j \left(\vartheta \right)-\tau _j\left(\vartheta
_0\right)<\delta   $. Then 
\begin{align*}
I_j\left(\vartheta \right)&\geq \int_{ \tau_j}^{\tau _j\left(\vartheta \right)} {S
  _{j}\left(t-\tau _j\right)} ^2{\rm d}t+\int_{ \tau_j\left(\vartheta
  \right)}^{\tau _j+\delta } \left[{S _{j}\left(t-\tau
    _j\left(\vartheta\right)\right)}- {S _{j}\left(t-\tau _j\right)}
  \right]^2{\rm d}t\\ 
&=\int_{0}^{\tau _j\left(\vartheta \right)- \tau_j} {S
  _{j}\left(s\right)} ^2{\rm d}s+ \int_{0} ^{\tau _j-\tau _j\left(\vartheta \right)+\delta }\left[S _{j}(s)- S _{j}\left(s-\Delta \tau _j\right)
  \right]^2{\rm d}s.
\end{align*}
where $\Delta\left( \tau _j\right)=\tau _j -\tau
    _j\left(\vartheta\right)$.
Further (below $\lambda _m=\min_{t\in {\cal T}_j}\lambda _j\left(t\right)>0$)
\begin{align*}
\int_{0} ^{\tau _j\left(\vartheta \right)-\tau_j} \lambda  _{j}\left(s\right)
^2\left(\frac{s}{\delta }\right)^{2\kappa }{\rm 
  d}s\geq \frac{\lambda _m^2}{\delta ^{2\kappa }}\int_{0} ^{\tau
  _j\left(\vartheta \right)-\tau_j}s^{2\kappa }{\rm 
  d}s=\frac{\lambda _m^2\left|\tau _j\left(\vartheta \right)-\tau_j
  \right|^{2\kappa +1} }{\delta ^{2\kappa }\left(2\kappa +1\right)}.
\end{align*}
Recall that 
\begin{align*}
\tau _j\left(\vartheta \right)-\tau_j=\langle m_j,\vartheta -\vartheta
_0\rangle +O\left(h^2\right)=\langle m_j,e\rangle \left\|\vartheta -\vartheta
_0\right\|+O\left(h^2\right),
\end{align*}
where the unit vector $e=\left(\vartheta -\vartheta
_0\right)\left\|\vartheta -\vartheta
_0\right\|^{-1} $. Therefore
\begin{align*}
&\int_{0} ^{\tau _j\left(\vartheta \right)-\tau_j} \lambda  _{j}\left(s\right)
^2\left(\frac{s}{\delta }\right)^{2\kappa }{\rm 
  d}s\\
&\qquad \qquad \geq \frac{\lambda _m^2\left|\langle m_j,e\rangle 
  \right|^{2\kappa +1} }{\delta ^{2\kappa }\left(2\kappa
    +1\right)}\left\|\vartheta -\vartheta 
_0\right\|^{2\kappa +1}\left(1+o\left(\left\|\vartheta -\vartheta
_0\right\|\right)\right)
\end{align*}
and we can take such $h$ that
\begin{align*}
\int_{0} ^{\tau _j\left(\vartheta \right)-\tau_j} \lambda  _{j}\left(s\right)
^2\left(\frac{s}{\delta }\right)^{2\kappa }{\rm 
  d}s\geq \frac{\lambda _m^2\left|\langle m_j,e\rangle 
  \right|^{2\kappa +1} }{2\delta ^{2\kappa }\left(2\kappa +1\right)}\left\|\vartheta -\vartheta
_0\right\|^{2\kappa +1}.
\end{align*}
For the second integral we have ($\delta _*=\tau _j-\tau _j\left(\vartheta
\right)+\delta >0$) 
\begin{align*}
&\int_{ 0}^{\delta_* } \left[{S
    _{j}\left(s\right)}- {S _{j}\left(s-\Delta \left(\tau _j\right)\right)} 
  \right]^2{\rm d}s\\
&\qquad =\frac{1}{\delta ^{2\kappa }}\int_{ 0}^{\delta_* } \left[{\lambda 
    _{j}\left(s\right)}s^\kappa - {\lambda  _{j}\left(s-\Delta \left(\tau _j\right)\right)\left|{s-\Delta \left(\tau _j\right)}\right|^\kappa} 
  \right]^2{\rm d}s\\
&\qquad \geq \frac{\lambda 
    _{m}^2}{\delta ^{2\kappa }}\int_{ 0}^{\delta_* } \left[s^\kappa -
    \left|s-\Delta \left(\tau _j\right)\right|^\kappa
  \right]^2{\rm d}s-C\left\|\vartheta -\vartheta
_0\right\|^2\\
&\qquad \geq\frac{\lambda 
    _{m}^2}{\delta ^{2\kappa }}\int_{ 0}^{\frac{\delta_* }{\Delta (\tau _j)}} \left[v^\kappa -
    \left|v-1\right|^\kappa
  \right]^2{\rm d}v\Delta \left(\tau _j\right)^{2\kappa +1} -C\left\|\vartheta -\vartheta
_0\right\|^2\\
&\qquad \geq\frac{\lambda 
    _{m}^2}{\delta ^{2\kappa }}\int_{ 0}^{\frac{\delta_* }{ch}} \left[v^\kappa -
    \left|v-1\right|^\kappa
  \right]^2{\rm d}v\left|\langle m_j,e\rangle 
  \right|^{2\kappa +1}\left\|\vartheta -\vartheta
_0\right\|^{2\kappa +1}\\
&\qquad \qquad  -C\left\|\vartheta -\vartheta
_0\right\|^2,
\end{align*}
where we used the relation $\lambda  _{j}\left(s-\Delta \left(\tau
_j\right)\right)=\lambda  _{j}\left(s\right)+O\left(\Delta \left(\tau
_j\right)\right) $ and set $s=v \Delta \left(\tau _j\right)$.

These estimates from below of the integral allow us to write
\begin{align*}
\sum_{j=1}^{k}I_j\left(\vartheta \right)\geq \gamma \sum_{j=1}^{k}\left|\langle m_j,e\rangle 
  \right|^{2\kappa +1}\left\|\vartheta -\vartheta
_0\right\|^{2\kappa +1}-C\left\|\vartheta -\vartheta
_0\right\|^2.
\end{align*}
As $k\geq 3$ we have
\begin{align*}
Q\left(e\right)= \sum_{j=1}^{k}\left|\langle m_j,e\rangle 
  \right|^{2\kappa +1},\qquad    \inf_{\left\|e\right\|=1}Q\left(e\right)= q_1>0.
\end{align*}
Indeed, if $q_1=0$, then there exists a vector $e_*$ such that
$Q\left(e_*\right)=0 $ and this vector is orthogonal to all
$m_j,j=1,\ldots,k$. Of course, this is impossible. Therefore we can take such
sufficiently small $h$ that for $\vartheta \in \AA\cap \BB$ we obtain the
estimate
\begin{align}
\label{ng}
\sum_{j=1}^{k}\int_{0}^{T}\left[{S 
  _{j}\left(t-\tau _j\left(\vartheta\right)\right)}- S 
  _{j}\left(t-\tau _j\left(\vartheta_0\right)\right) \right]^2{\rm d}t\geq \gamma _1
\left\|\vartheta -\vartheta _0\right\|^{2\kappa +1}
\end{align}
with some positive $\gamma _1$. For the other values of $\vartheta \in \AA$ we
have the similar estimates.

Let us consider these integrals for the  values $\vartheta \in
\AA^c$.  According to \eqref{bel} we have to study
 the function
\begin{align*}
g\left(h\right)=\inf_{\vartheta _0\in\Theta
}\inf_{\left\|\vartheta -\vartheta _0 
  \right\|>h }\sum_{j=1}^{k}\int_{0}^{T}\left[{S 
  _{j}\left(t-\tau _j\left(\vartheta\right)\right)}- S 
  _{j}\left(t-\tau _j\left(\vartheta_0\right)\right) \right]^2{\rm d}t ,
\end{align*}
 and show that $g\left(h\right)>0$.

   Suppose that $g\left(h \right)=0$, then this implies
  that there exists at least one point $\vartheta^*\in \Theta $ such that
  $\left\|\vartheta ^*-\vartheta _0\right\|\geq h$ and for all
  $j=1,\ldots,k$ we have
\begin{align*}
\int_{0}^{T}{\left[S  
  _{j}\left(t-\tau _j\left(\vartheta^* \right)\right)- S
  _{j}\left(t-\tau _j\left(\vartheta_0 \right) \right)\right]^2}{\rm d}t=0.
\end{align*}
Let $\tau _j\left(\vartheta ^*\right)>\tau _j$. Then   for all $t\in
\left[\tau _j,\tau _j\left(\vartheta ^*\right) \right]$ we have
\begin{align*}
\lambda   _{j}\left(t-\tau _j \right)\left|{t-\tau _j }\right|^\kappa = 0
\end{align*}
and for $t\in
\left[\tau _j\left(\vartheta ^*\right),\tau _j+\delta  \right] $
\begin{align*}
\lambda _{j}\left(t-\tau _j\left(\vartheta^* \right)\right)\left|{t-\tau
  _j\left(\vartheta ^*\right)}\right|^\kappa = \lambda  _{j}\left(t-\tau
_j \right)\left|{t-\tau _j }\right|^\kappa.
\end{align*}
Of course we can have these two equalities if and only if $\tau
_j\left(\vartheta ^*\right)= \tau _j\left(\vartheta_0\right)$ for all
$j=1,\ldots,k$. Recall that $\lambda _j\left(t\right)$ are strictly positive
functions. From the geometry of the model it follows that it is impossible to
have two different points such that the distances from these points and $k\geq 3$
detectors coincide.

Therefore for $\vartheta \in \AA^c$
\begin{align}
\label{k}
&\sum_{j=1}^{k}\int_{0}^{T}\left[{S
  _{j}(t-\tau _j\left(\vartheta\right))}-S
  _{j}(t-\tau_j \left(\vartheta_0\right)) \right]^2{\rm d}t \geq g\left(h
\right)\nonumber\\
&\qquad \qquad \geq  \frac{g\left(h
\right)\left\|\vartheta -\vartheta _0    \right\|^{2\kappa +1}}{ D^{2\kappa
    +1}}\geq \gamma _2 \left\|\vartheta -\vartheta _0    \right\|^{2\kappa +1},
\end{align}
where $D=\sup_{\vartheta _1,\vartheta _2\in\Theta
}\left\|\vartheta_1-\vartheta _2\right\|$.

From the estimates \eqref{ng} and \eqref{k} it follows that if we put
$\vartheta =\vartheta _0+\nu \varphi _nu$, then 
\begin{align*}
\sum_{j=1}^{k}\int_{0}^{T}\left[\sqrt{\lambda
  _{j,n}\left(\theta_u,t\right)}- \sqrt{\lambda
  _{j,n}\left(\vartheta_0,t\right)} \right]^2{\rm d}t  & \geq \gamma
n\left\|\vartheta -\vartheta _0    \right\|^{2\kappa +1}\\
&=\gamma \nu ^{2\kappa
  +1} \left\|u   \right\|^{2\kappa +1} .
\end{align*}
This estimate and \eqref{32} prove \eqref{pi}.

\bigskip

The properties of the likelihood ratio field $Z_n\left(\cdot \right)$
established in the lemmas \ref{L1}-\ref{L3} are sufficient  conditions for the
Theorem 1.10.2 in \cite{IH81}. Therefore the  Theorem \ref{T1} is proved.

\section{Discussion}

There are several problems which naturally arise for this model of
observations. Note that the properties of the MLE $\hat\vartheta _n$ can
be studied too. This requires a special modification of Lemma 2 to verify the
tightness of the corresponding family of measures. 

 In this work we supposed
that the source starts  emission at the instant $t=0$. It is
interesting to consider the more general statement with unknown start of the
emission. The limit distributions of the MLE and BE are unknown and it will be
interesting to have some pictures obtained by numerical simulations for the
densities of these vectors. The numerical simulations can provide us the values
 $\Ex_{\vartheta _0}\|\hat\zeta \|^2$   and
$\Ex_{\vartheta _0}\|\tilde\zeta \|^2$   of the  limit variances
of these estimators. 

Note that it is possible to construct a consistent estimator of $\vartheta _0$
in two steps as it was proposed in \cite{CK18}. First we estimate $k$ moments
$\tau =\left(\tau _1,\ldots,\tau _k\right)$ of
arriving signals in detectors, say,  $\tilde \tau_ {1,n},\ldots,\tilde \tau_
{k,n}$. Recall that
\begin{align*}
\tilde\xi _{j,n}=n^\frac{1}{2\kappa +1}\left( \tilde\tau _{j,n}-\tau _j\right)\Longrightarrow \tilde\xi _j,\quad j=1,\ldots,k
\end{align*}
where $\tilde\xi _j$ are independent random variables. Hence we have
\begin{align*}
\nu ^2\tilde\tau _{j,n}=\nu ^2\tau _j^2+2\varphi _n\tau _j{\tilde\xi
  _{j,n}}+\varphi _n^2{\tilde\xi _{j,n}}^2 =\rho _j^2+2\nu \rho _j\tilde\xi
  _{j,n}\varphi _n +O\left(\varphi _n^2\right).
\end{align*}
 Then we write the equations  
\begin{align*}
\left(x_j-x_0^*\right)^2+\left(y_j-y_0^*\right)^2=\nu ^2\tilde\tau _{j,n}=\rho _j^2+2\nu \rho _j\tilde\xi
  _{j,n}\varphi _n +O\left(\varphi _n^2\right),\quad j=1,\ldots,k,
\end{align*}
 and obtain the least squares estimator $\vartheta _n^*$, which is
consistent and has the same rate of convergence as the BE $\tilde\vartheta _n$
studied in this work. See details in the Section 3, \cite{CK18}.

{\bf Acknowledgment.} This work was done under partial financial support of
the grant of RSF number 14-49-00079.


\begin{thebibliography}{99}

%
\bibitem{CK18} Chernoyarov, O.V. and Kutoyants, Yu.A. (2018)
  Source localization on the plane. Smooth case. To be submitted.


\bibitem{D03} Dachian, S. (2003) Estimation of cusp location by Poisson
  observations. {\sl Statist. Inference Stoch. Processes.} 6, 1, 1-14.

\bibitem{D11} Dachian, S. (2011) Estimation of the location of a 0-type or
  $\infty $-type singularity by Poisson observations. {\sl A Journal of
  Theoretical and Applied Statistics,} 45, 5, 509-523.

\bibitem{FKT18} Farinetto, C.,  Kutoyants, Yu.A. and Top, A. (2018)  
  Source localization on the plane. Change point case. To be submitted.
%


%
\bibitem{IH81} Ibragimov, I. A. and Khasminskii, R. Z. (1981)  {\it Statistical
Estimation. Asymptotic Theory.} New York: Springer.
%
\bibitem{Ka91} Karr, A. F. (1991) 
{\it Point Processes and their Statistical Inference.} New York: Marcel Dekker.
%



\bibitem{Kn10} Knoll, G.F.  (2010) {\it Radiation Detection and
  Measurement.} New York: Wiley. 


\bibitem{Kut79} Kutoyants Yu. A. (1979)  Parameter estimation of intensity of
inhomogeneous Poisson processes. {\it Problems of Control and Information
Theory,}  8, 137-149.

\bibitem{Kut84} Kutoyants Yu. A. (1984)  {\it Parameter Estimation for Stochastic
  Processes.} Heldermann: Berlin.  
%
\bibitem{Kut98} Kutoyants, Yu. A. (1998)  {\it Statistical Inference for
  Spatial Poisson Processes}. New York: Springer.

%

\bibitem{LS01} Liptser, R.S., Shiryayev, A.N. (2001) \textsl{Statistics of
  Random Processes, I, II.} (2nd edition), Springer, New York.
 
\bibitem{XL13} Luo, X. (2013) {\it GPS Stochastic Modelling.} New York: Springer.


%
\bibitem{SM91} Snyder, D.R. and Miller, M.I. (1991) 
{\it Random point processes in time and space.}
New York: Springer.
%

\bibitem{St10} Streit, R.L. (2010)  {\it Poisson Point Processes: Imaging,
  Tracking, and Sensing.}  Boston: Springer.

%

\end{thebibliography}
\end{document}